\documentclass[12pt,thmsa]{article}
\usepackage{amsfonts}
\usepackage[dvips]{graphics}

\def\mylabel#1{\label{#1}}
\newtheorem{theorem}{Theorem}[section]
\newtheorem{lemma}[theorem]{Lemma}
\newtheorem{corollary}[theorem]{Corollary}
\newtheorem{proposition}[theorem]{Proposition}
\newtheorem{example}[theorem]{Example}
\newtheorem{remark}[theorem]{Remark}

\newtheorem{hypothesis}[theorem]{Hypothesis}

\newtheorem{definition}[theorem]{Definition}
\def\bit{\begin{itemize}}
\def\eit{\end{itemize}}
\reversemarginpar   
\def\bc{\begin{center}}
\def\ec{\end{center}}
\def\bthm{\begin{theorem}}
\def\ethm{\end{theorem}}
\def\bcor{\begin{corollary}}
\def\ecor{\end{corollary}}
\def\bprop{\begin{proposition}}
\def\eprop{\end{proposition}}
\def\blem{\begin{lemma}}
\def\elem{\end{lemma}}
\def\bex{\begin{example} }
\def\eex{\end{example} }
\def\brem{\begin{remark}}
\def\erem{\end{remark}}
\def\prf{\noindent{\bf Proof~: }}
\def\bdes{\begin{description}}
\def\edes{\end{description}}
\def\ita{\item[(a)]}
\def\itb{\item[(b)]}
\def\itc{\item[(c)]}

\def\iti{\item[(i)]}
\def\itii{\item[(ii)]}
\def\itiii{\item[(iii)]}
\def\itiv{\item[(iv)]}
\def\itv{\item[(v)]}

\def\beq{\begin{equation}}
\def\eeq{\end{equation}}
\def\ben{\begin{enumerate}}
\def\een{\end{enumerate}}
\def\beqar{\begin{eqnarray}}
\def\eeqar{\end{eqnarray}}
\def\beqarr{\begin{eqnarray*}}
\def\eeqarr{\end{eqnarray*}}

\def\RR{{\mathbb R}}  

 \def\cE{\mathcal{E}} \def\cF{\mathcal{F}}  
    
  \def\cL{\mathcal{L}}

\def\qed{\hspace{.1in}{\bf QED}}

\def\P{{\mathsf P}} 

\def\E{{\mathsf E}} 

\def\ZZ{{\mathbb Z}}       

\def\NN{{\mathbb N}}       

\def\one{{\bf 1}}

\def\rar{\rightarrow}

\def\eps{\epsilon}

\def\part{\partial}

\def\d#1dt{\frac{d#1}{dt}}    

\def\hrar{\hookrightarrow}

\begin{document}
\title{A  class of non homogeneous self interacting random processes with
applications to Learning in Games and Vertex-Reinforced Random Walks\thanks{We acknowledge financial support from the Swiss National Science Foundation grant 200020-112316.}}
\author{{\bf Michel Bena\"{\i}m}\\ 
Universit\'e de Neuch\^atel, Suisse \and
{\bf Olivier Raimond}\\ 
Universit\'e Paris Sud, France} 
\maketitle
\bibliographystyle{apalike}
\begin{abstract}
Using an approximation by a set-valued dynamical system,
this paper studies a class of non Markovian and non homogeneous
stochastic processes on a finite state space. It provides an unified
approach to {\em simulated annealing} type processes. It permits to
study new  models of {\em vertex reinforced random walks} and new
models of learning in games including {\em Markovian fictitious play}.

\paragraph{Keywords:} Stochastic approximation, Processes with reinforcement, Differential Inclusions, Learning in Games, Simulated Annealing.
\end{abstract}

\section{Introduction}
\label{intro}
Let  $E$ be  a finite set  called {\em the state space,} $\mathsf{M} =
 \mathsf{M}(E)$   the set of Markov matrices over $E,$ and  $\Sigma$
 a compact convex subset of an Euclidean space called  {\em the
 observation space.}
The set $\Sigma$ will be equipped with the distance induced by the
 Euclidean norm $\|\cdot\|$ on the observation space.
 Let $(\Omega,\cF, \P)$ be a probability space equipped with an increasing sequence of sub $\sigma$-fields $\{\cF_n, n \in \NN\}: \: \cF_n \subset \cF_{n+1} \subset \cF.$

Our main object of interest is  a  discrete time  random process $(X,M,V) = ((X_n, M_n,V_n))$ defined on $(\Omega,{\cal F}, \P)$  taking values in $E \times {\mathsf{M}}(E) \times \Sigma$ such that:
 \bdes
\iti $(X,M,V)$ is {\em adapted} (to $\{\cF_n, n \in \NN\}$), meaning that $(X_n,M_n,V_n)$ is ${\cal F}_n$-measurable for each $n.$
\itii For all $y \in E$
\beq
\label{eq:markov} \P(X_{n+1} = y |{\cal F}_n) = M_n(X_n,y).
\eeq
\edes
We refer to $X_n$ (respectively $V_n$) as the {\em state} (respectively, the {\em observation}) variable  at time $n;$ and to the sequence $(M_n)$ as the {\em strategy.} 
We let $$v_n = \frac{1}{n} \sum_{i = 1}^n V_i$$ denote  the empirical average up to time $n$ of the sequence of observations.

A well studied situation is when 
\beq
\label{eq:control}
M_n = K(v_n)
\eeq
where $K$ maps continuously probability vectors to irreducible Markov matrices and
$$V_{n+1} = H(X_{n+1},v_n)$$ for some map $H : E \times \Sigma \mapsto \Sigma.$ In such a case  $(X_n)$ is called a ``Markov chain controlled'' by $(v_n)$ and the behavior of $(v_n)$  can be analyzed through the ODE 
\beq
\label{eq:ode0}
\dot{v} = -v + \sum_x \pi(v)(x)H(x,v)
\eeq where $\pi(v)$ is the invariant probability of $K(v).$
This approach to controlled Markov chains goes back to the work of M\'etivier and Priouret (1987) (see also the books  Benveniste, M\'etivier and Priouret (1990), Duflo (1996))   strongly influenced by   the pioneered works of Ljung (1977),  Kushner and Clark (1978) on the ODE's method. It has been used in Bena\"{\i}m (1997) for analyzing certain vertex reinforced random walks on finite graphs.

The main purpose of this paper is to investigate the long term behavior of $(v_n)$ under 
less stringent assumptions than (\ref{eq:control}). In particular we are interested in situations where: \bdes
\ita  $M_n$ may depend 
on other (non-observable or hidden)  variables than $v_n$ and;
\itb  The  closure of $\{M_n\::n \geq 0\}$ may contain degenerate (i.e~ non irreducible) Markov matrices.
\edes
Situation (a)  typically occurs in game theory where players may have only partial information on the actions played by their opponents, and  (b) is motivated by stochastic optimization algorithms.

Relying on a recent paper by Bena\"{\i}m, Hofbauer and Sorin (2005)  it will be shown that under certain assumptions (involving estimates on the log-Sobolev and spectral gap constants of $(M_n)$) the asymptotic behavior can be described in term of a certain set-valued deterministic dynamical system that generalizes the ODE (\ref{eq:ode0}). Applications to non-homogeneous Markov chains, vertex reinforced random walks and learning processes in game theory will be given.

\subsection*{Outline of contents} The organization of the paper is as follows.
Section \ref{sec:notation} states the notation, hypotheses and  the
main result. Our main assumption (Hypothesis \ref{hyp:main}) is somewhat abstract and more tractable conditions  (expressed in term of spectral gaps and log-Sobolev constants) are given in section \ref{sec:mainhyp}.
Section  \ref{sec:appli} is devoted to examples and applications. The proof of the main result is postponed to section \ref{sec:proof}.
\section{Notation, hypotheses and main results}
\mylabel{sec:notation}
A {\em probability vector} (or measure) over $E$ is a map $\mu : E \to \RR^+$ such that $\sum_x \mu(x) = 1,$ and a  {\em Markov matrix}  is a map $M : E \times E \to \RR^+,$ such that
$$\forall x \in E, \: \sum_{y} M(x,y) = 1.$$
We let
$\Delta = \Delta(E)$ denote the space of probability vectors over $E$ and 
$\mathsf{M} = \mathsf{M}(E)$ denote  the set of Markov matrices on $E.$ 

Given a function $f : E \to \RR$ and  $\mu \in \Delta$ we use the notation
$$\mu f = \sum_x \mu(x) f(x).$$
A Markov matrix $M$ on $E$ acts on functions $f$ and measures $\mu$ according to the formulas
$$Mf(x) = \sum_y M(x,y) f(y),$$
$$\mu M(y) = \sum_x \mu(x) M(x,y).$$
We let $M^n$ denote the Markov matrix obtained by matrix multiplication. Equivalently $M^n f = M (M^{n-1} f)$ for $n \geq 1,$ with the convention that $M^0f = f.$

Points $x, y \in E$ are said to be {\em related} if there exist $i, j \geq 0$ (depending on $x$ and $y$) such that $M^i(x,y) > 0$ and $M^j(y,x) > 0.$ An equivalence class for this relation  is called a {\em recurrent class.}
The Markov matrix $M$ on $E$ is said {\em indecomposable} if it has a unique recurrent class (possibly periodic) and is said {\em irreducible} if this recurrent class is $E.$ 

By standard results, indecomposability of $M$ implies that $M$ possesses a  unique  {\em invariant probability measure} $\pi$ characterized by  the relation $\pi M = \pi.$ Moreover, the {\em generator} $L = -I + M$ has kernel $\RR 1$ and its restriction to  $\{f \::  \pi f = 0\}$ is an isomorphism. It then follows that $-L$ admits a  pseudo ``inverse'' $Q$ characterized by
$$Q\one = 0,$$
and
$$Q (I-M)  = (I-M) Q = I -\Pi;$$ 
where $\Pi \in \mathsf{M}$ denote the matrix defined by  $\Pi(x,y) = \pi(y).$
To shorten notation we also call $Q$ the pseudo inverse of $M.$

Given a vector $f$ and a matrix $N,$   we set
$|f| = \max |f(x)|$ and $|N| = \max_{x,y} |N(x,y)|.$

\medskip
Our main assumption is the following:
\begin{hypothesis}
\mylabel{hyp:main}
The matrices $(M_n)$ are  indecomposable and their pseudo inverses $(Q_n)$ and invariant probabilities $(\pi_n)$ satisfy almost surely
\bdes
\iti
$$\lim_{n \to \infty} \frac{|Q_n|^2 \log(n)}{n} = 0,$$
\itii $$\lim_{n \to \infty} |Q_{n+1}-Q_n| = 0,$$
\itiii $$\lim_{n \to \infty} |\pi_{n+1}-\pi_n| = 0.$$
\edes
\end{hypothesis} 

The verification of  hypothesis \ref{hyp:main} is the subject of section \ref{sec:mainhyp} where  sufficient and more tractable conditions will be detailed. 

\medskip
Let $\hat{V}_n : E \to \Sigma$ be an $\cF_n$-measurable map defined by
$$\hat{V}_n(x) = \frac{\E(V_{n+1}\one_{X_{n+1} = x} | \cF_n)}{M_n(X_n,x)}$$
for $M_n(X_n,x) \neq 0.$
In addition to hypothesis \ref{hyp:main} we assume that
\begin{hypothesis}
\mylabel{hyp:main2}
$$\lim_{n \to \infty} M_{n+1} Q_{n+1} (\hat{V}_{n+1} - \hat{V}_n) = 0$$
almost surely.
\end{hypothesis}

\brem
\mylabel{remhyp2}
 {\rm Here are some sufficient conditions ensuring hypothesis \ref{hyp:main2}.
\bdes
\iti 
Assume that $x \mapsto \hat{V}_{n+1}(x) - \hat{V}_n(x)$ is a constant map. 
Then hypothesis \ref{hyp:main2} holds since $Q_{n} 1 = 0$. This will be used in section \ref{sec:appli}. 

\itii More generally, let $T \Sigma$ be the affine hull of $\Sigma$ (the smallest affine space containing $\Sigma$). Assume that 
for all $n \in \NN$ there exists a vector $A_n \in T \Sigma$ and a map  $B_n : E \to  T \Sigma$ such that
\bdes
\ita For all $x \in E$, $\hat{V}_{n+1}(x) - \hat{V}_n(x) = A_n + B_n(x)$
\itb $\limsup_{n \to \infty} |B_n|
 \sqrt{\frac{n}{\log(n)}} < \infty$, almost surely. 
\edes
Then $|M_{n+1}Q_{n+1} ( (\hat{V}_{n+1} - \hat{V}_n))| = |M_{n+1}Q_{n+1} B_n|
\leq |Q_{n+1}||B_n| \to 0$ almost surely by hypothesis \ref{hyp:main}.

\itiii Assume that $M_n(x,y) = \pi_n(y).$ Then $M_{n+1} Q_{n+1} = 0$ so that hypothesis \ref{hyp:main2} holds. \edes

} \erem
\subsection{Adapted set-valued dynamical systems}
The purpose of this section is to introduce certain  differential inclusions on $\Sigma$ that will prove to be useful for analyzing the long term behavior of $(v_n).$
Recall that we let $\pi_n$ denote the invariant probability of $M_n.$
Let
\beq
\label{eq:thetan}
\theta_n =  \pi_n \hat{V}_n = \sum_x \pi_n(x) \hat{V}_n(x).
\eeq
We let $C_n \subset \Sigma \times \Sigma$ denote the {\em topological support} of the law of $(v_n,\theta_n).$ That is  the smallest closed set $F \subset \Sigma \times \Sigma$ such that
$$\P((v_n,\theta_n) \in F) = 1.$$
Let $\mathsf{clos}{\{C_n\}}$ denote the set of all possible limit points $z = \lim z_{n_k}$ with $z_{n_k} \in C_{n_k}$ and $n_k \rar \infty.$ It is easily seen that $\mathsf{clos}{\{C_n\}}$ is a nonempty compact subset of $\Sigma \times \Sigma$. 

A nonempty set $G \subset \Sigma \times \Sigma$ is called a {\em graph} (or a bundle) over $\Sigma$, if the projection $$p : G \to \Sigma,$$
$$(u,v) \mapsto u$$ is onto.
A graph $G$ over $\Sigma$ defines a {\em set-valued function} mapping each point $u  \in \Sigma$ to a set $G(u) = \{v \in \Sigma \: : (u,v) \in G\}.$ 
\begin{definition}
\label{defC}
A set $C \subset \Sigma \times \Sigma$ is said to be {\em adapted} to  $\{(v_n,\theta_n)\}$ (or simply adapted) if
\bdes
\iti
$C$ is a closed graph over $\Sigma.$
\itii For all $u \in \Sigma, \, C(u)$ is a nonempty convex set.
\itiii  $\mathsf{clos}{\{C_n\}} \subset C.$
\edes
\end{definition}

\medskip
To an  adapted set $C$ we associate the  differential inclusion 
\beq
\label{eq:inclus}
\dot{v} \in -v + C(v).
\eeq
A solution to (\ref{eq:inclus}) is an absolutely continuous mapping $v:\RR \to \Sigma$ verifying $\dot{v}(t) + v(t) \in C(v(t))$ for almost every $t.$
A set $A \subset \Sigma$ is said to be  
 {\em  invariant}  if for all $x \in A$ there exists a solution
  ${\bf x}$ to (\ref{eq:inclus}) with ${\bf x }(0) = x
$ and  such that ${\bf x}(\RR) \subset A.$

Given a set $A \subset \Sigma$ and $(x, y) \in A^2$ we write
$x \hrar_A y$ if for every ${\varepsilon} > 0$ and $T > 0$ there exists an
  integer $n \in \NN,$ solutions ${\bf x}_1, \ldots {\bf x}_n$ to
(\ref{eq:inclus})
and real numbers
  $t_1, t_2, \ldots, t_n $ greater than $T$ such that
\bdes
\ita
${\bf x}_i([0,t_i]) \subset A,$
\itb $\|{\bf x}_i(t_{i}) - {\bf x}_{i+1}(0)\| \leq {\varepsilon} $ for all $i
= 1,
\ldots,n-1,$
\itc $\|{\bf x}_1(0) - x\| \leq {\varepsilon}$  and $\|{\bf x}_n(t_n) - y\|
\leq {\varepsilon}.$
\edes
\begin{definition}  A set $A \subset \Sigma$ is said to be {\em internally chain
transitive}  provided $A$ is compact  and
$x \hrar_A y$ for all $x, y \in A.$
\end{definition}
It is not hard to verify (see e.g~Bena\"{\i}m, Hofbauer and Sorin (2005) Lemma 3.5) that an {\em internally chain
transitive} set is invariant.

The {\em limit set} of $(v_n)$ is the set  $L = L((v_n))$ consisting of all points $p = \lim v_{n_k}$ for some sequence $n_k \to \infty.$ The next theorem \ref{th:main} is the main result of the paper. Its proof heavily  relies on Bena\"{\i}m, Hofbauer and Sorin (2005) and is given in section \ref{sec:proof}.
\bthm
\mylabel{th:main}
Assume that hypotheses  \ref{hyp:main} and  \ref{hyp:main2} hold. Let $C$ be an adapted graph. Then the limit set of $(v_n)$ is an internally chain transitive set for the differential inclusion
$$\dot{v} \in -v + C(v).$$
\ethm
\subsection{Background : How to use Theorem \ref{th:main}}
The notion of ``internally chain transitive set''  was introduced by Bena\"{\i}m and Hirsch (1996) in order to analyze the long term behavior of certain perturbations of flows and has been  recently extended to multivalued dynamical systems by Bena\"{\i}m, Hofbauer and Sorin (2005). We refer the reader to this paper for more details, examples and properties. For convenience this section briefly reviews   a few  useful properties of internally chain transitive sets. 

The differential inclusion (\ref{eq:inclus}) induces a
  set-valued dynamical system $\{\Phi_t\}_{t \in \RR} $ defined by
$$\Phi_t(x) = \{{\bf x}(t) : \: {\bf x} \mbox{ is a solution to
(\ref{eq:inclus}) with } 
  {\bf x} (0) = x \in \Sigma\}.$$
 
A non empty compact
set  $A$ is  an {\em attracting} set if there exists
 a neighborhood $U$  of $A$ and a function ${\bf t}$ from $(0,
{\varepsilon}_0)$ to $\RR^+$  with ${\varepsilon}_0> 0$
such that
$$
\Phi_t(U) \subset A^{{\varepsilon}}
$$
 for all
 $\varepsilon <  \varepsilon_0$ and
$t \geq {\bf t} (\varepsilon )$, where $A^{{\varepsilon}}$
 stands for the ${\varepsilon}-$neighborhood of $A.$  
If additionally $A$ is invariant, then $A$ is  an  {\it attractor}.

Given an attracting set (resp. attractor) $A,$ its {\it basin of attraction} is the set
  $$B(A) = \{x \in \Sigma: \: \exists t \geq 0, \, \Phi_t(x) \in U\}.$$ When $B(A) = \Sigma,$  $A$   is a   {\it globally}
 attracting set (resp. a global attractor).

Given a closed invariant set $S,$ the induced  dynamical
system $\Phi^S$ on $S$ is defined by
$$
\Phi^S_t(x) = \{{\bf x}(t) :  {\bf x} \mbox{ is a solution to }
(\ref{eq:inclus}) \mbox{ with } {\bf x} (0) = x   \mbox{ and } {\bf x}(\RR) \subset S\}.
$$
An invariant  set $S$ is  {\it attractor free} if  there exists  no proper subset $A$ of $S$  which is  an attractor for $\Phi^S$.\\

Throughout the remainder of this section we let $L$ denote an internally chain transitive set (for instance the limit set  $L = L(v_n)$).
Properties of  $L$ will then be obtained
through the next  result  (Bena\"{\i}m, Hofbauer and Sorin (2005), Lemma 3.5, Proposition 3.20 and
 Theorem 3.23):
\bprop
\mylabel{propertyB}
\bdes
\iti The set $L$ is   non--empty, compact,  invariant  and   attractor free.
\itii  If $A$ is an attracting set with $B(A) \cap L  \not=
 \emptyset,$  then $L \subset  A$.
\edes
\eprop

Some useful properties of attracting sets
or attractors are the two following (Bena\"{\i}m, Hofbauer and Sorin (2005), Propositions 3.25 and 3.27).
\bprop
\mylabel{propertyC}  
 Let $\Lambda \subset \Sigma$ be  compact   with  a bounded 
open
neighborhood $U$ and $V : \overline{U} \to [0, \infty[$. Assume the following conditions:
\bdes
 \iti $\Phi_t(U) \subset U$ for all $t \geq 0$,
\itii  $V^{-1}(0) = \Lambda,$
\itiii $V$ is continuous and   for all $x \in U \setminus \Lambda,
 y \in \Phi_t(x)$ and $t > 0,$  $V(y)< V(x).$
\edes
Then  $ \Lambda$ contains an attractor whose basin 
contains $U.$
\eprop

The map $V$ introduced in this proposition is called a {\em strong Lyapounov function} associated to $\Lambda.$

Let now $\Lambda$ be a subset of $\Sigma$ and $U \subset \Sigma$ an open neighborhood of $\Lambda.$ A continuous function $V : U \to \RR$ is called a 
{\em Lyapunov function} for $\Lambda \subset \Sigma$ if 
$V(y) < V(x)$ for all $x \in U \setminus \Lambda$, $y \in \Phi_t(x)$, $t > 0;$ 
and $V(y) \leq V(x)$ for all $x \in \Lambda$,  $y \in \Phi_t(x)$ and 
$t \geq 0.$ 
\bprop[Lyapounov]
\mylabel{propertyD}
Suppose $V : U \to \RR $ is a Lyapunov function for $\Lambda$ and $L\subset U.$ 
Assume that $V(\Lambda)$ has an empty interior.
Then $L \subset  \Lambda$ and the restriction of $V$ to $L$ is constant.
\eprop
\section{Verification of hypothesis \ref{hyp:main}}
\mylabel{sec:mainhyp}
This section is devoted to the verification of Hypothesis \ref{hyp:main}. The results given here  will be used in section \ref{sec:appli} to analyze specific situations.  

\subsection{Estimates based on compactness}
Let $\mathsf{M}_{ind}(E)$ denote the open set of indecomposable Markov matrices.
\bprop 
\mylabel{th:compact2}
Suppose that the sequence $(M_n)$ lies in a compact subset of
$\mathsf{M}_{ind}(E)$ and verifies $\lim_{n \to \infty} (M_{n+1}-M_n) = 0.$  Then hypothesis \ref{hyp:main} holds.
\eprop
This proposition is a direct consequence of the  next lemma.

\blem
\mylabel{th:compact}
Let $T\mathsf{M}(E)$ be the space of matrices $K = K(x,y)$ such that $\sum_y K(x,y) = 0.$
The map $\mathsf{Q} : \mathsf{M}_{ind}(E) \to T\mathsf{M}(E)$ which associates to $M$ its pseudo inverse  and the map $\mathsf{\Pi} : \mathsf{M}_{ind}(E) \to \Delta$ which associates to $M$ its invariant measure are  smooth maps. 
\elem
\prf Set $M \in \mathsf{M}_{ind}(E)$. The invariant probability of $M$, $\Pi(M)$, is solution to $\phi(M,\pi)=0$ where
$\phi :  \mathsf{M}_{ind}(E) \times  \Delta \to  T\Delta,$ is the smooth map defined by
$$\phi (M,\mu) =   \mu (I-M),$$
with $T\Delta=\{\mu:E\to\RR~:\sum_x \mu(x)=0\}$. For all $\nu \in  T\Delta$, $$\frac{\partial \phi}{\partial \mu}(M,\mu).\nu =  \nu (I-M).$$
Hence, by uniqueness of the invariant probability measure, $\frac{\partial \phi}{\partial \mu}(M,\mu)$ has kernel $\{0\}$ and the fact that $\Pi$ is smooth follows from the implicit function theorem.
 
We denote by $\hat{\Pi}(M) \in \mathsf{M}(E)$ the matrix defined by  $\hat{\Pi}(M)(x,y) = \Pi(M)(y).$ The pseudo inverse of $M$ is solution to $\psi(M,Q) = 0$ where  $\psi :  \mathsf{M}_{ind}(E) \times  T\mathsf{M}(E) \to  T\mathsf{M}(E),$ is the smooth map defined by
$$\psi (M,Q) =   Q (I-M) - (I- \hat{\Pi}(M)).$$ 
For all $A \in  T\mathsf{M}(E)$ $$\frac{\partial \psi}{\partial Q}(M,Q).A =  A (I-M).$$
Hence, 
by uniqueness of the invariant probability measure, $\frac{\partial \psi}{\partial Q}(M,Q)$ has kernel $\{0\}$ and the fact that $Q$  depends smoothly on $M$ follows from the implicit function theorem.  \qed

\medskip
Let $K$ be a continuous mapping from $\Gamma$ a
compact set into $\mathsf{M}(E)$ such that $K(w)$ is indecomposable
for all $w\in\Gamma$. Assume $(w_n)$ is a sequence of $\Gamma$-valued
random variables such that $M_n=K(w_n)$. 
If in addition $\lim_{n\to\infty}(M_{n+1}-M_n)=0$, then
proposition \ref{th:compact2} applies.

\subsection{Estimates based on log-Sobolev and spectral gap constants}
Propositions \ref{th:hypmain} and   \ref{th:ls} below can be used to verify hypothesis \ref{hyp:main}   when the sequence  $(M_n)$ is not bounded away from $\mathsf{M}_{ind}(E).$ The strategy is then to verify assertions $(ii)$ and $(iii)$ of proposition \ref{th:hypmain} and to  use the estimates given by proposition \ref{th:ls} to verify assertion $(i).$

\bprop
\mylabel{th:hypmain}
Suppose that the matrices $(M_n)$ are  indecomposable and that their pseudo inverse $(Q_n)$ and invariant probabilities $(\pi_n)$ satisfy amost surely
\bdes
\iti
$$\lim_{n \to \infty} \frac{|Q_n|^2 \log(n)}{n} = 0,$$
\itii $$\limsup_{n \to \infty} |M_{n+1} - M_n| \frac{n}{\log(n)} < \infty$$
\itiii $$\limsup_{n \to \infty} |\pi_{n+1} - \pi_n| \sqrt{\frac{n}{\log(n)}} < \infty.$$
\edes
Then hypothesis \ref{hyp:main} holds.
\eprop  
\prf The proof amounts to show that hypothesis \ref{hyp:main} $(ii)$ holds.
Set $L_n = M_n - I$ and $\Pi_n=\hat{\Pi}(M_n)$. Using the characterization of $Q_n$ one has
 $$Q_{n+1} (L_{n+1} - L_n) + (Q_{n+1} - Q_n) L_n = \Pi_{n+1} - \Pi_n.$$ Hence, 
$$Q_{n+1} (L_{n+1} - L_n) Q_n  + (Q_{n+1} - Q_n) L_n Q_n = (\Pi_{n+1} - \Pi_n) Q_n.$$
That is (using $Q_n\Pi_n=Q_n\Pi_{n+1}=0$ and $L_nQ_n=\Pi_n-I$)
$$Q_{n+1} (M_{n+1} - M_n) Q_n  + (Q_{n} - Q_{n+1})  = (\Pi_{n+1} - \Pi_n) Q_n.$$ 
Therefore
$$|Q_n - Q_{n+1}| \leq   c (|Q_{n+1}||Q_n||M_{n+1} - M_n| + |\pi_{n+1} - \pi_n||Q_n|),$$
for some constant $c > 0$
and conditions
 $(i), (ii), (iii)$ imply hypothesis \ref{hyp:main} (ii).
\qed

\medskip
Let  $\mathsf{M}_{irr}(E)$ denote the open set of irreducible Markov matrices.
Let $M \in \mathsf{M}_{irr}(E)$ with invariant probability $\pi$ and let  $f : E \to \RR.$ 
The {\em variance}, {\em entropy} and  {\em energy} of $f$
are respectively  defined as
$$var(f) = \pi (f^2) - (\pi f)^2$$
$$\cL(f) = \sum_x f(x)^2 \log\left(\frac{f(x)^2}{\pi f^2}\right) \pi(x)$$
$$\cE(f) = \frac{1}{2} \sum_{x,y} (f(y)-f(x))^2 M(x,y) \pi(x).$$
The {\em spectral gap} and {\em log-Sobolev} constants of $M$ are then defined to be
$$\lambda = \min \left\{\frac{\cE(f)}{var(f)} \: : var(f) \neq 0\right\}$$
$$\alpha =  \min \left\{\frac{\cE(f)}{\cL(f)} \: : \cL(f) \neq 0\right\}.$$
The following estimates follows from the quantitative results for finite Markov chains  as given in Saloff-Coste (1997) theorems. 

\bprop
\mylabel{th:ls}
 Let $M \in \mathsf{M}_{irr}(E)$ with invariant probability $\pi$ log-Sobolev constant $\alpha$ and spectral gap $\lambda.$ For all $(x,y) \in E$ 
the following estimates hold:
\bdes
\iti $$|Q(x,y)| \leq \sqrt{\frac{\pi(y)}{\pi(x)}}\frac{1}{\lambda}$$
\itii $$|Q(x,y)| \leq \frac{1}{\alpha} \log_{+} \left(\log\left(\frac{1}{\pi(x)}\right)\right) + \frac{e}{\lambda}$$
where $\log_{+}(t) = \max(0,\log(t)).$
\edes
In particular
$$|Q| \leq \frac{1}{\alpha}  \left [  \log_{+} \left(\log\left(\frac{1}{\pi_*}\right)\right) + \frac{e}{2}\right ]$$ and 
$$|Q| \leq \frac{1}{\lambda} \left [ \log_{+} \left(\log\left(\frac{1}{\pi_*}\right)\right) \frac{\log((1-\pi_*)/\pi_*)}{1-2\pi_*}) + e \right ].$$
\eprop
\prf
Let $L = -I + M$ and let $\{P_t\}$ be the continuous time semi-group  $P_t = e^{t L}.$ Then $Q$ can be written as
$$Q(x,y) =  \int_0^{\infty} (P_t(x,y) - \pi(y))dt.$$
The first assertion then easily follows from the estimate $$|P_t(x,y) -\pi(y)| \leq 
 \sqrt{\frac{\pi(y)}{\pi(x)}} e^{-\lambda t}$$ whose proof can be found in Saloff-Coste (1997, Corollary 2.1.5).

We now pass to the second assertion. If $\pi(x) \geq e^{-2}$ the inequality to be proved follows from inequality $(i).$  Hence we assume that  $\pi(x) < e^{-2},$ and we  follow the line of the proof of Theorem 2.2.5 in Saloff-Coste (1997).  For $q \geq 1,$ we let $||.||_{q}$ denotes the norm in $l^q(\pi).$ We let $P_t^*$ denote the adjoint of $P_t$ in $l^2(\pi),$ and  $p_t(x,y) = p_t^*(y,x) = P_t(x,y)/\pi(y).$ Let $g_x$ denote the function given by $g_x(y) = 0$ for $x \neq y$ and $g_x(x) = 1/\pi(x).$ Then
$$|P_t(x,y) - \pi(y)| \leq 
||p_t(x,.) -1||_2 = ||(P_t^* - \pi)g_x||_2$$
 Therefore 
$$|P_{t+s}(x,y) - \pi(y)| \leq ||p_{t+s}(x,.) -1||_2  \leq ||P_t^*-\pi||_{2 \to 2} ||P_s^* g_x||_2$$
$$\leq e^{-\lambda t} ||P_s^*||_{k \to 2} ||g_x||_k$$ for any $k \geq 1.$
where we have used the fact that $||P_t^*-\pi||_{2 \to 2} \leq e^{-\lambda t}.$
Let $q$ be the H\"older conjugate of $k.$ Then  $||P_s^*||_{k \to 2} =  ||P_s||_{2 \to q}.$ Now choose $q(s) = 1+e^{2\alpha s}.$ By hypercontractivity (see Theorem 2.2.4 in Saloff-Coste (1997)), $||P_s||_{2 \to q(s)} \leq 1$ so that
$$|P_{t+s}(x,y) - \pi(y)| \leq e^{-\lambda t}  \pi(x)^{-1/q(s)}.$$
Hence
$$|Q(x,y)| \leq  2 s +  \frac{1}{\lambda}  \pi(x)^{-1/q(s)}.$$ For $s = \frac{1}{2\alpha} \log_{+} (\log(\frac{1}{\pi(x)}))$ this gives the desired result.

The uniform bounds on $|Q|$ follow from the rough estimates 
$$\frac{1-2\pi_*}{\log((1-\pi_*)/\pi_*)} \lambda \leq \alpha \leq \lambda/2$$
given in Saloff-Coste (1997, Lemma 2.2.2 and Corollary 2.2.10) \qed
\section{Some applications}
\mylabel{sec:appli}
In sections \ref{sec:markov2} and \ref{sec:vrrw2},
we are interested in the long term behavior of the {\em empirical occupation measure of the process}. We then let $\Sigma = \Delta, V_n = \delta_{X_n}$ and 
$$v_n = \frac{1}{n} \sum_{i = 1}^n \delta_{X_i}.$$
Hence,  $\hat{V}_n(x) = \delta_x$ and  $\theta_n = \pi_n.$ 
\subsection{Markov chains}
\mylabel{sec:markov2}
Let $(M_n)$ be a deterministic (or $\cF_0$ measurable) sequence of Markov matrices over $E.$ 
A {\em non homogeneous} Markov chain with transition matrices $(M_n)$ is an adapted process $(X_n)$ on $E$ verifying (\ref{eq:markov}). 


 
\bprop
\mylabel{th:cormark1}
Let $L((\pi_n)) \subset \Delta $ denote the limit set of $(\pi_n)$ and let $\mathsf{conv}[L((\pi_n))]$ denote its convex hull. 
Suppose that hypothesis \ref{hyp:main} holds.  Then
$L((v_n)) \subset \mathsf{conv}[L((\pi_n))]$  with probability one.
\eprop
\prf The set
$C = \Delta \times \mathsf{conv}[L((\pi_n))]$ is adapted to $(v_n,\pi_n)$. The induced differential equation $\dot{v} \in  - v + \mathsf{conv}[L((\pi_n))]$ has a unique global attractor $\mathsf{conv}[L((\pi_n))].$ Hence, by Theorem \ref{th:main} and Proposition \ref{propertyB}, (ii), $L((v_n)) \subset \mathsf{conv}[L((\pi_n))].$
\qed
\bcor Suppose that the sequence $(M_n)$ lies in a compact subset of
$\mathsf{M}_{ind}(E)$ and verifies $M_{n+1}-M_n \to 0.$  Then conclusion of proposition \ref{th:cormark1} holds. \ecor
\prf Follows from proposition \ref{th:cormark1} and proposition \ref{th:compact2}. \qed
\bcor Assume that $M_n \to M \in \mathsf{M}_{ind}(E).$ Then $v_n \to \pi$ the invariant probability of $M.$ \ecor
\subsubsection*{Markov chains with rare transitions}
Among the well studied chains that motivate our analysis are the {\em chains with rare transitions.}

Let $M_0$ be an irreducible  Markov matrix over $E,$ reversible with respect to a reference probability $\pi_0.$ 
That is $$\pi_0(x) M_0(x,y) = \pi_0(y) M_0(y,x).$$
We sometimes call such an $M_0,$ an  {\em exploration matrix} since it
 provides a way to explore the state space.

Let $W : E \times E \to \RR,$ be a map and $(\beta_n)$ a sequence of positive numbers.
Set 
\beq
\label{eq:simann1}
M_n(x,y) = M(\beta_n,x,y)
\eeq
where
$$
M(\beta,x,y) =   \left  \{
\begin{array}{lll} M_0(x,y)\psi[ \exp (-\beta W(x,y))] & \mbox{if } x \neq y, \\ & \\
1-\sum_{y \neq x} M(\beta,x,y) & \mbox{if } x =  y,
\end{array} \right.
$$
and
\beq
\label{eq:psi}
\psi(u) = \min (1, u)
\eeq
 or $$\psi(u) =  \frac{u}{1+u}.$$
In particular, let $U : E \to \RR$
be  a map, and let 
\beq
\label{eq:simann2} 
W(x,y) = U(y) - U(x),  
\eeq
then $(M_n)$ are the transition matrices of the so-called {\em Metropolis-Hasting}  ($\beta_n = \beta$) or {\em simulated annealing} ($\beta_n \to \infty$) algorithm (Hajek (1982), Holley and Stroock (1988), Miclo (1992)).

Consider the Markov chain with rare transitions (\ref{eq:simann1})
where $W$ is given by (\ref{eq:simann2}).
For $x,y \in E$ a path $\gamma$ from $x$ to $y$ is a sequence of points
$x_0 = x, x_1, \ldots x_n = y$ such that $M_0(x_i,x_{i+1}) > 0.$
We let $\Gamma_{x,y}$ denote the set of all paths from $x$ to $y.$ The {\em elevation} from $x$ to $y$ is defined as 
$$\mathsf{Elev}(x,y) = 
\min \{ \max \{U(z) \: : z \in \gamma\} \: : \gamma \in \Gamma_{x,y}\}$$ and the {\em energy barrier} as 
\beq
\label{eq:barr}
U^{\#} = \max \{\mathsf{Elev}(x,y) - U(x) - U(y) + \mathsf{min}U \: : x \in E, y \in E \}
\eeq

\bprop 
\mylabel{th:recuit}
Consider the Markov chain with rare transitions (\ref{eq:simann1})
with $W$ given by (\ref{eq:simann2}).
Assume that $\beta_n = \beta(n)$ where $\beta : \RR^+ \to \RR^+$
is differentiable and verify
 $$ 0 \leq \dot{\beta}(t) \leq \frac{A}{t}$$
for some $A < 1 / 2U^{\#}.$ 
Then $v_n \to \pi$
where
$$\pi(x) \propto \pi_0(x) \one_{\mathsf{Argmin} U}(x).$$   
\eprop
\prf Our first goal is to verify hypothesis \ref{hyp:main}.
Let $\lambda(\beta)$ denote the spectral gap of $M(\beta,\cdot,\cdot).$
It follows from Theorem 2.1 in Holley and Stroock (1988) that
\beq
\label{eq:holley}
\lim_{\beta \to \infty} \frac{\log(\lambda(\beta))}{\beta} = - U^{\#}.
\eeq
The invariant probability measure of $M(\beta,\cdot,\cdot)$ is the {\em Gibbs measure}
\beq
\label{eq:gibbs}
\pi_{\beta}(x) \propto \exp (-\beta U(x))\pi_0(x).
\eeq
Since $\beta_n \leq \beta_1 + A \log(n),$ by application of the last inequality of Proposition \ref{th:ls},
one gets that hypothesis \ref{hyp:main} (i) holds.

For $x \neq y$
 $$\frac{\partial M(\beta,x,y)}{\partial \beta} = - M_0(x,y)W(x,y)
\psi'(\exp(-\beta W(x,y)) \exp(-\beta W(x,y)).$$
Using the fact that   $|\psi'(t)t| \leq 1,$ one gets that
 $$ \left|\frac{\partial M(\beta,x,y)}{\partial \beta}\right| \leq c$$
for some $c > 0.$
Hence by the mean value theorem
$$|M_{n+1} - M_n| \leq c |\beta_{n+1} - \beta_n| \leq (A c)/n.$$
This proves assertion $(ii)$ of proposition \ref{th:hypmain}.
The proof of assertion $(iii)$ is similar since 
 $$\left|\frac{\partial \pi_\beta(x)}{\partial \beta}\right| = |\pi_{\beta}(x) (U(x) - \sum_y \pi_{\beta}(y) U(y))|   \leq 2||U||.$$ 
This concludes the verification of  hypothesis \ref{hyp:main}.

Here $\pi_n(x)  \propto \exp (-\beta_n U(x)) \pi_0(x)$ so that 
$\pi_n \to \pi.$ The result  follows from 
Proposition \ref{th:cormark1}.
\qed

\brem {\em For general $W,$ it is always possible to define a  {\em quasipotential} 
$U$ (defined in term of $W$ and $M_0$) and an energy barrier  $U^{\#}$ (in general not given by (\ref{eq:barr})) such that both equations
(\ref{eq:holley}) and (\ref{eq:gibbs})  hold.  We refer the reader to Miclo (1992) for more details and proofs. With this quasi-potential and barrier  
Proposition \ref{th:recuit} holds.}
\erem
\subsection{Vertex reinforced random walks}
\mylabel{sec:vrrw2}
Vertex-reinforced random walks (VRRW) were first introduced by Pemantle (1988, 1992).

Suppose $\cF_n = \sigma(X_1,\ldots, X_n).$ A general VRRW  on $E$
is  defined
by 
$$M_n(x,y) = K_n(x,y,v_n)$$
where for each integer $n$  and $v \in \Delta$, 
$K_n(\cdot,\cdot,v)$ is a deterministic Markov matrix over $E,$ which specifies the rule of the reinforcement.

The following result was proved in Bena\"{\i}m (1997).
\bprop
\mylabel{th:vrrw1}
Assume that there exists a $[0,1]$-valued sequence $\epsilon_n$ converging to $0$ at infinity such that $K_n(x,y,v) = K(x,y,\epsilon_n,v)$, that the map $(\epsilon,v) \mapsto K(\cdot,\cdot,\epsilon,v)$ is continuous on $[0,1]\times \Delta$ and that $K(\cdot,\cdot,\epsilon,v)$ is indecomposable for each $(\epsilon,v) \in [0,1]\times \Delta$. Let  $\pi(v)$ denote the invariant measure of $K(x,y,0,v)$. Then the limit set of $(v_n)$ is almost surely an internally chain transitive set of the differential equation
\beq
\label{eq:odevrrw} \dot{v} = -v + \pi(v).
\eeq
\eprop 
\prf This follows from  Proposition \ref{th:compact2} and Theorem \ref{th:main}. \qed
\subsubsection*{Linear reinforcement}
The original VRRW  as defined by Pemantle (1988, 1992) corresponds to a {\em linear reinforcement:} 
$$M_n(x,y) \propto U(x,y) \left[1 + \sum_{i = 1}^n \one_{X_i = y}\right],$$
where $U$ is a  matrix with nonnegative entries. 

We will here assume that $U$ has positive entries.
Then, for each $n$, $M_n$ is irreducible.
With the notation of the previous paragraph,
\beq
\label{eq:knpem1}
M_n(x,y) = K(x,y,1/n,v_n),
\eeq
where for $(\epsilon,v)\in [0,1]\times\Delta$,
\beq
\label{eq:knpem1b}
K(x,y,\epsilon,v) \propto
U(x,y) \left[\epsilon + v(y)\right].
\eeq
The mapping $(\epsilon,v)\mapsto K(\cdot,\cdot,\epsilon,v)$ is continuous on $[0,1]\times \Delta$. 

On a finite graph, this process was first analyzed by Pemantle (1992) for symmetric positive matrices  ($U(x,y) = U(y,x) > 0$)  and later by Bena\"{\i}m (1997) for general positive matrices using proposition \ref{th:vrrw1}.
As an example of what can be proved is the following result first due to Pemantle (1992)
\bprop\label{prop:lyap}
Suppose $U(x,y) = U(y,x) > 0.$ Then the limit set of $(v_n)$ is a compact connected subset of the critical set of the map
$$v \mapsto U(v,v) = \sum_{x,y} U(x,y)v(x) v(y).$$ 
\eprop
\prf 
This follows from the fact that $v \mapsto U(v,v)$ is a strict lyapounov function of (\ref{eq:odevrrw})  whose critical points are the zeroes of (\ref{eq:odevrrw}).
\qed 

\medskip
When the matrix $U$ has zero entries, $K(x,y,0,v)$ may no longer be indecomposable for some $v \in  \partial\Delta$ and proposition \ref{th:vrrw1} cannot be applied. This makes the analysis of VRRW with linear reinforcement much more difficult.  Beautiful results  on $\ZZ$ and $\ZZ^d$ have been obtained by Pemantle and Volkov (1999), Volkov (2001) and Tarres (2004). We refer the reader to Pemantle (2007) for a survey and further references. 

\medskip

\subsubsection*{Non homogeneous linear reinforcement}
Let $(a_n)$ be a positive sequence and denote $r_n=\sum_{i=1}^n
a_i$. We will assume that $\lim_{n\to\infty}\frac{r_{n+1}}{r_n}=1$.
Consider the VRRW corresponding to:
$$M_n(x,y) \propto U(x,y) \left[1 + \sum_{i = 1}^n a_i \one_{X_i = y}\right],$$
where $U$ is a  matrix with positive entries.
Equivalently, $M_n(x,y)=K(x,y,\eps_n,w_n)$ with
\beq
\label{eq:knpem2}
K(x,y,\eps,w) \propto U(x,y) \left[\eps + w(y)\right],
\eeq
$\eps_n=1/r_n$ and $w_n=\frac{1}{r_n}\sum_{i=1}^n
a_i\delta_{X_i}$. Using proposition \ref{th:compact2}, it is not hard
to check that hypothesis \ref{hyp:main} and hypothesis
\ref{hyp:main2} (with $V_i=\delta_{X_i}$) are satisfied, so that
theorem \ref{th:main} applies.

\medskip
Since
$\delta_{X_i}=v_i+(i-1)(v_i-v_{i-1})$, using the convention $r_0=v_0=0$,
\beqarr
w_n
&=& \frac{1}{r_n}\sum_{i=1}^n (r_i-r_{i-1}) v_i + \frac{1}{r_n}\sum_{i=1}^n (i-1)(v_i-v_{i-1}) a_i \\
&=& v_n + \frac{1}{r_n}\sum_{i=1}^{n-1} r_i(v_i-v_{i+1})  + \frac{1}{r_n}\sum_{i=1}^{n-1} ia_{i+1}(v_{i+1}-v_i) \\
&=& v_n -\frac{1}{r_n}\sum_{i=1}^n(r_i-ia_{i+1})(v_{i+1}-v_i).
\eeqarr
Since $|v_{i+1}-v_i|\leq 2/i$,
$$|w_n-v_n|\leq
\frac{2}{r_n}\sum_{i=1}^n\left|\frac{r_i}{i}-a_{i+1}\right|.$$

Consider now the two following classes of sequences $(a_i)$:
\bdes
\iti $a_i=a(i)$ where $a$ is a nondecreasing continuous function such that
for all positive $s\in ]0,1]$, $\lim_{t\to\infty}\frac{a(ts)}{a(t)}=1.$
\itii $a_i=a(i)$ where $a$ is a decreasing continuous function such that
for all positive $s\in ]0,1]$,
$\lim_{t\to\infty}\frac{a(ts)}{a(t)}=1$,
and there exists $b:[0,1]\to \RR^+$ measurable such that $\int_0^1
b(s)ds<\infty$ and for all $(s,t)\in ]0,1]\times\RR^+$,
$$0\leq \frac{a(ts)}{a(t)}-1\leq b(s).$$
\edes
For example $a_i=(\log(i+1))^\alpha$ satisfies (i) for $\alpha\geq 0$ and (ii) for $\alpha<0$.

\blem Assume (i) or (ii) holds, then $\lim_{n\to\infty}|w_n-v_n|=0$. \elem
\prf Note that it suffices to prove that $\frac{r_i}{i}-a_{i+1}= o(a_i)$. 
Assume first(i) holds. Then
\beqarr
0&\leq& a_{i+1}-\frac{r_i}{i} \\
&\leq& a_i\left( \frac{a_{i+1}}{a_i}-1 + \int_0^1 \left(1-\frac{a(is)}{a(i)}\right)ds\right)\\
&=& o(a_i).
\eeqarr

Assume now (ii) holds. Then 
\beqarr
0&\leq& \frac{r_i}{i} -a_{i+1} \\
&\leq& a_i\left( 1-\frac{a_{i+1}}{a_i} + \int_0^1 \left(\frac{a(is)}{a(i)}-1\right)ds\right)\\
&=& o(a_i)
\eeqarr
by dominated convergence theorem. \qed
\medskip

Let $\pi(\epsilon,v)$ denote the invariant probability of
$K(x,y,\epsilon,v)$ and $\pi(v)=\pi(0,v)$. The map $(\epsilon,v)\mapsto \pi(\epsilon,v)$ is uniformly continuous. Then the previous lemma implies that when (i) or (ii) holds, since $\pi_n=\pi(\epsilon_n,w_n)$, $\lim_{n\to\infty}|\pi_n-\pi(v_n)|=0$. This last property with theorem \ref{th:main} implies the
\bthm Assume that (i) or (ii) holds, then
the limit set of $(v_n)$ is almost surely an internally chain transitive set of the differential equation
\beq
\dot{v} = -v + \pi(v).
\eeq
\ethm
Note that proposition \ref{prop:lyap} also holds for sequences $(a_i)$ satisfying (i) or (ii).

\subsubsection*{Exponential reinforcement}
Let $U : E \times E \to \RR$ be a map. For  $x \in E$ and $v \in \Delta,$
set
$$U(x,v) = \sum_{y \in E} U(x,y)v(y),$$ 
$$ 
W(x,y,v) = U(y,v) - U(x,v),
$$
$$
K(\beta,x,y,v) =   \left  \{
\begin{array}{lll} M_0(x,y)\psi[ \exp (-\beta W(x,y,v))] & \mbox{if } x \neq y, \\ & \\
1-\sum_{y \neq x} K(\beta,x,y,v) & \mbox{if } x =  y,
\end{array} \right.
$$
and
\beq
\label{eq:vrrw2} 
K_n(x,y,v) = K(\beta_n,x,y,v),
\eeq
Here $M_0$ is an exploration matrix, $(\beta_n)_n$ is a positive sequence and $\psi$ is given by (\ref{eq:psi}).
When $\beta_n = \beta$, such a VRRW can be seen as a discrete time version of the self-interacting diffusions on compact manifolds that have been thoroughly analyzed by 
Bena\"{\i}m, Ledoux and Raimond (2002),  Bena\"{\i}m and Raimond (2003, 2005, 2006).
When $\beta_n=A\log(n)$, the VRRW can be seen as a discrete time version of the self-interacting diffusions on compact manifolds studied by Raimond (2006).

\medskip
Let $U^{\#}(\cdot,y)$ be the energy barrier as defined by equation (\ref{eq:barr}) of the map
$x \mapsto   U(x,y)$ 

\bthm 
\mylabel{th:vrrw2}
Consider the VRRW with exponential reinforcement defined by (\ref{eq:vrrw2}).
Assume that $\beta_n = \beta(n)$ where $\beta : \RR^+ \to \RR^+$
is differentiable and verify
 $$ 0 \leq \dot{\beta}(t) \leq \frac{A}{t}$$
for some $A < 1 /  2  \max\{U^{\#}(\cdot,y) \: : y \in E\}.$
Let $$C(v) = \Delta(\mathsf{Argmin} U(\cdot, v))$$ 
denote the set of probabilities supported by $\mathsf{Argmin} U(\cdot, v).$
Then the limit set of $(v_n)$ is an internally chain transitive set of
$$\dot{v}\in -v+C(v).$$
\ethm
\prf This is an application of Theorem \ref{th:main}.
The verification of hypothesis \ref{hyp:main} is similar to the one given in proposition
 \ref{th:recuit}. Details are left to reader.

It is easily seen that $C$ is a closed-valued set with convex values.
For $v \in \Delta,$ let $$\pi_n[v](x) \propto \pi_0(x) \exp(-\beta_n U(x,v))$$
and $$\pi[v](x) \propto \pi_0(x) \one_{\mathsf{Argmin} U(\cdot, v)}(x).$$
The invariant probability of $K_n$ is $\pi_n[v_n]$  and
$$\lim_{n \to \infty} \pi_n[v](x) = \pi[v](x).$$  
This proves that $C$ is adapted to $(v_n,\pi_n[v_n])$ and the result follows from Theorem \ref{th:main}.
\qed 

\bcor[symmetric interaction]
\mylabel{th:vrrw3}
Assume that hypotheses of  Theorem \ref{th:vrrw2} hold and 
  assume furthermore that  $U$  is {\em symmetric}  (i.e~$U(x,y) = U(y,x)$).
Then $(v_n)$ converges almost surely to a connected component of the set 
$$\{v \in \Delta: \: v \in C(v)\}.$$
\ecor
\prf 
For $u,v \in \Delta$ set $$U(u,v) = \sum_{x,y} U(x,y)u(x)v(y)$$ and 
let $$H(v) = \frac{1}{2}U(v,v)$$
We claim that $H$ is a lyapouvov function of the differential inclusion
(\ref{eq:inclus}).
Let $t  \mapsto v(t)$ be a solution to (\ref{eq:inclus}) then, for almost all $t \geq 0$ 
\begin{eqnarray*}
\frac{d}{dt} H(v(t)) &=& \frac{1}{2}[U(\dot{v}(t),v(t)) + U(v(t),\dot{v}(t)]
 = U(\dot{v}(t),v(t))\\
&=& U(\dot{v}(t)+v(t),v(t)) -U(v(t),v(t))\\
&=& \min_{x} U(x,v(t)) -U(v(t),v(t)),
\end{eqnarray*}
where we have used the symmetry of $U,$  the fact that $\dot{v} + v \in C(v)$ and the definition of $C(v).$
Since $t \mapsto H(v(t))$ is locally Lipchitz, it is nondecreasing. If now   $t \mapsto H(v(t))$ is constant over a time interval, then $v(t) \in C(v(t))$ over this time interval. This proves that $H$ is a Lyapounov function for $\Lambda = \{v \in \Delta \: : v \in C(v)\}.$  The result now follows from Proposition \ref{propertyD} (compare to Bena\"{\i}m, Hofbauer and Sorin (2005), Theorem 5.5) provided we show that $H(\Lambda)$ has empty interior. 

Let $v \in \Lambda \cap \mathsf{int}(\Delta)$. Since the mapping $x\mapsto U(x,v)$ is constant, for all for all $w\in \Delta$, $U(w,v)=U(v,v)$. Therefore $H(v) = U(w,v)$ for all $w \in \Delta$.  It follows that $H$ restricted to $\Lambda \cap \mathsf{int}(\Delta)$ is a constant map. The same reasonning applies to prove that $H$ restricted to each face of $\Delta$ is a constant map. We thus have proved that $H(\Lambda)$ takes finitely many values. \qed

\brem {\em Corollary \ref{th:vrrw3} still holds true under the weaker assumption that the map $v \mapsto U(x,v)$ is smooth and convex in $v.$}
\erem

\bcor Assume that $U$ is symmetric and nonnegative and that $$Ker(U) \cap T\Delta = \{0\}.$$ Then $\{v \in \Delta: \: v \in C(v)\}$ reduces to a singleton $v^*$ and $(v_n)$ converges almost surely to $v^*.$ \ecor
\prf Let $v \in C(v), w \in \Delta$ and $h = w-v.$ Since $v \in C(v), U(v,h) \geq 0.$  Thus $U(w,w) - U(v,v) = 2 U(v,h) + U(h,h) \geq 0,$  proving that $v$ is a global minimum of $v \mapsto U(v,v).$ Since $U(h,h) > 0$ for $h = w -v \neq 0,$ such a global minimum is unique. \qed

\subsection{Games}
Consider a  two-players game. We let
$E_1$  (respectively $E_2$) denote  the finite  {\em set of actions}
 available to player $1$  (respectively  player $2$)  and $$U = (U^1, U^2) : E_1 \times E_2 \to \RR \times \RR$$ 
denote the payoff function of the game. If player $1$ and player $2$ 
 choose respectively the actions $x \in E_1$ and $y \in E_2$, then player $1$ gets   $U^1(x,y)$ and player $2$ gets $U^2(x,y)$.

Let $((X_n,Y_n))$ denote the sequence of plays. In noncooperative game theory we assume that $((X_n,Y_n))$ is  adapted to some filtration $(\cF_n)$  and that at the beginning of round $n+1,$  players have no information on the action to be played by their opponents: for all $(x,y)\in E_1\times E_2$ and $n\in\NN$
$$\P(X_{n+1} = x, Y_{n+1} = y |\cF_n) = \P(X_{n+1} = x |\cF_n)\P(Y_{n+1} = y |\cF_n).$$
\subsubsection{Markovian fictitious play}
For $x \in E_1$ and $v^2 \in \Delta(E_2)$   set 
$$U^1(x,v^2) = \sum_{z \in F} U^1(x,z)v^2(z).$$
Let 
$$v^2_n = \frac{1}{n} \sum_{i = 1}^n \delta_{Y_i}.$$ A well studied strategy  known as ``fictitious play''  consists for player $1$ to play at time $n+1$ an action  maximizing $U^1( \cdot, v^2_n)$,
that is
\beq
\label{eq:fictmax}
X_{n+1} \in \mathsf{Argmax} U^1(\cdot,v^2_n ).
\eeq

This strategy relies on the idea that in absence of information on the next move of his opponent, player 1 assumes that he (the opponent)  will play accordingly to the past empirical distribution of his  moves.
While fictitious play was originally proposed in 1951 by Brown as an algorithm to compute Nash equilibria it has been recently rediscover  as a ``learning model'' (Fudenberg and Kreps (1993); Fudenberg and Levine (1998)) and has been extensively studied (Monderer and Shapley (1996); Bena\"{\i}m and Hirsch (1999); Hofbauer and Sandholm (2002); Bena\"{\i}m, Hofbauer and Sorin (2005, 2006), see also Pemantle (2007) for an  overview  and further references).

Fictitious plays requires to solve the maximization problem (\ref{eq:fictmax}) at each stage of the game. 
If the cardinal of  $E_1$ is too large (or if players have computational limits) such a computation may be problematic. An alternative strategy proposed first in Bena\"{\i}m, Hofbauer and Sorin (2006), based on pairwise comparison of payoffs, is as follows: The strategy of player 1 is such that
$\P(X_{n+1}=y|\cF_n)=M_n(X_n,y)$
with $M_n$ the Markov matrix defined by
\beq
\label{eq:fict1} 
M_n(x,y) =   \left  \{
\begin{array}{lll} M_0(x,y)\psi[ \exp (-\beta_n W_n(x,y))] & \mbox{if } x \neq y, \\ & \\
1-\sum_{y \neq x} M_n(x,y) & \mbox{if } x =  y,
\end{array} \right.
\eeq
where 
$$W_n(x,y) = U^1(x,v^2_n) - U^1(y,v^2_n),$$
$M_0$ is an exploration matrix, $\psi$ is given by (\ref{eq:psi}) and $\beta_n$ is an increasing positive sequence.
Such a strategy will be called a {\em Markovian fictitious play strategy}.



Adopting the view point of player 1, 
 we choose, as an observation space,
  $$\Sigma = \Delta(E_1) \times \Delta(E_2)$$
and as an observation variable
$$V_n = (\delta_{X_n}, \delta_{Y_n}).$$ Hence $(v_n)$ is the {\em empirical frequency of the actions} played up to time $n,$ and 
$$\hat{V}_n(x) = (\delta_x, \nu_n),$$
where $\nu_n = \E(\delta_{Y_{n+1}}|\cF_n).$

\medskip
We let $U^{1,\#}(y)$ denote the energy barrier, as defined by (\ref{eq:barr}), of the map $x \mapsto U^1(x,y).$

\bthm
\mylabel{th:fic1}
Assume that player 1 plays a Markovian fictitious play strategy as given by (\ref{eq:fict1}). Assume that 
$\beta_n = \beta(n)$ where $\beta$ is differentiable, $\lim_{t \rar \infty}  \beta(t) = \infty$ and verify
$$0 \leq \dot{\beta}(t) \leq \frac{A}{t}$$ for some $A <1/ 2 \max\{U^{1,\#}(y) \: : y \in E_2\}.$

For $v = (v^1,v^2) \in \Delta(E_1) \times \Delta(E_2)$  let 
 $$C_1(v^2) = \Delta (\mathsf{Argmax} U^1( \cdot, v^2))$$
and 
$$C(v) = C_1(v^2) \times \Delta(E_2)$$
Then the limit set of $(v_n)$ is an internally chain transitive set of
$$\dot{v}\in -v + C(v).$$
\ethm
\prf This is still an application of Theorem \ref{th:main}. The verification of hypothesis \ref{hyp:main} is similar to the one given in proposition \ref{th:recuit}. 
Let $$\pi_n[v^2](x) \propto \pi_0(x) \exp (\beta_n U^1(x,v^2))$$
and $$\pi[v^2](x) \propto \pi_0(x) \one_{\mathsf{Argmax}(U^1(\cdot,v^2))}(x).$$
Then, the invariant probability of $M_n$ is $\pi_n = \pi_n[v^2_n]$ and  
$\theta_n = \pi_n \hat{V}_n =  (\pi_n, \nu_n)$ with $\nu_n = \E(\delta_{Y_{n+1}}|\cF_n).$ Since $\pi_n[v^2] \to \pi[v^2] \in C^1(v^2)$ it follows that $C$ is an adapted graph. \qed 

\medskip
Much more can be said under the assumption that {\bf both} players adopt a Markovian fictitious play strategy: $\P(X_{n+1}=y|\cF_n)=M^1_n(X_n,y)$ and $\P(Y_{n+1}=y|\cF_n)=M^2_n(Y_n,y)$,
with $M^1_n$ and $M^2_n$ the Markov matrices defined by (with $i\in\{1,2\}$)
\beq
\label{eq:fict2} 
M^i_n(x,y) =   \left  \{
\begin{array}{lll} M_0(x,y)\psi[ \exp (-\beta^i_n W^i_n(x,y))] & \mbox{if } x \neq y, \\ & \\
1-\sum_{y \neq x} M^i_n(x,y) & \mbox{if } x =  y,
\end{array} \right.
\eeq
where 
\begin{eqnarray*}
W^1_n(x,y) &=& U^1(x,v^2_n) - U^1(y,v^2_n),\\
W^2_n(x,y) &=& U^2(v^1_n,x) - U^2(v^1_n,y),
\end{eqnarray*}
$M^i_0$ is an exploration matrix, $\psi$ is given by (\ref{eq:psi}) and $\beta_n^i$ is an increasing positive sequence.

Let $\mathsf{Conv}(U)$ denote the convex hull in $\RR^2$ of the set $\{U(x,y) : x \in E_1, y \in E_2\}$ of all possible payoffs. We now choose 
$$\Sigma = \Delta(E_1) \times \Delta(E_2) \times
 \mathsf{Conv}(U)$$  as an  observation space,
and 
$$V_n = (\delta_{X_n}, \delta_{Y_n}, U(X_n,Y_n))$$ as the observation variable. Hence
$$\hat{V}_n(x,y) = (\delta_x, \delta_y,  U(x,y)).$$
\bthm 
\mylabel{th:fic2}
Assume that both players adopt a Markovian fictitious play strategy.
Assume that for $i\in\{1,2\}$, 
$\beta^i_n = \beta^i(n)$ where $\beta^i$ is differentiable and verify
$$0 \leq \dot{\beta^i}(t) \leq \frac{A^i}{t}$$ for some $A^i <1/ 2 \max\{U^{i,\#}(y) \: : y \in E_{3-i}\}.$
 
For $v = (v^1,v^2,u) \in \Delta(E_1) \times \Delta(E_2) \times
 \mathsf{Conv}(U),$  let 
$$C(v) = \{\{ (\alpha,\beta, \gamma) \in \Sigma : \alpha \in C_1(v^2), \beta \in C_2(v^1), \gamma = U(\alpha,\beta)\}$$ 
where   $C_1(v^2)$ is like in Theorem \ref{th:fic1} and $C_2(v^1)$ is analogously defined for player 2.
Then the limit set of $(v_n)$ is an internally chain transitive set of $$\dot{v}\in -v+C(v).$$
\ethm
\prf
Let $(M_n^i)$  denote the strategy of Player i. Let $\pi_n^i, \lambda_n^i$ be the invariant measure and spectral gap of $M_n^i.$ 
On the  state space $E_1 \times E_2$
 the strategy of the pair of players is  $M_n = M_n^1 \otimes M_n^2$
which invariant measure is $\pi_n = \pi_n^1\otimes \pi_n^2$ and spectral gap
 $\lambda_n = \min(\lambda_n^1,\lambda_n^2).$ Thus hypothesis \ref{hyp:main} holds for $(M_n).$ The rest of the proof is similar to the proof of Theorem \ref{th:fic1} and is left to the reader. \qed
 
\bcor[zero sum games] Suppose that $U^2 = -U^1.$ Then under the assumption of Theorem \ref{th:fic2}, $(v^1_n,v_n^2)$ converges almost surely to the set of {\em Nash equilibria}
$$\{(v_1,v_2) \: :v_1 \in C_1(v^2), v_2 \in C_2(v^1) \},$$
and $(U^1(X_n,Y_n))$ converges almost surely to the {\em value of the game} 
$$u^* = \max_{v^1 \in \Delta(E_1)} \min_{v^2 \in \Delta(E_2)} U^1(v^1,v^2) = \min_{v^1 \in \Delta(E_1)}
 \max_{v^2 \in \Delta(E_2)} U^1 (v^1,v^2).$$
\ecor
\prf This follows from theorem \ref{th:main}, proposition \ref{propertyB} (ii) and the fact that the set $\{(v_1,v_2,u) \: :v_1 \in C_1(v^2), v_2 \in C_2(v^1),u = u^* \}$ is a global attractor of the differential inclusion, as proved in full generality by Bena\"{\i}m, Hofbauer and Sorin (2005). \qed

\bcor[Potential games]  Suppose that $U^2 = U^1.$ Then under the assumption of Theorem \ref{th:fic2}, $(v^1_n,v_n^2)$ converges almost surely to a connected subset  of the set of Nash equilibria
$$\{(v_1,v_2) \: :v_1 \in C_1(v^2), v_2 \in C_2(v^1) \}$$
on which $U^1$ is constant, 
and $(U^1(X_n,Y_n))$ converges almost surely towards this constant. 
\ecor
\prf Follows from theorem \ref{th:main}, proposition \ref{propertyD},  and the fact that $U^1 = U^2$ is a Lyapounov function of the differential inclusion. The proof of this later point is given in    (Bena\"{\i}m Hofbauer and Sorin, 2005, Theorem 5.5). It is similar to the proof  Corollary \ref{th:vrrw3}.\qed
\subsubsection{A remark on hypothesis \ref{hyp:main2}}
We give here a simple example showing the necessity of hypothesis \ref{hyp:main2}.

Consider the zero sum game where $E_1 = E_2 = \{0,1\}$, $U^1 = -U^2$ and
$$
U^1= \left [ \begin{array}{r r}
 0  & -1 \\ 
-1  &  0   
\end{array} \right ].
$$
Let $V_n = U^1(X_n,Y_n)$ be the payoff to player 1 at time $n.$
One has $$\hat{V}_n(x) =  U^1(x,1)\nu_n + U^1(x,0)(1-\nu_n)$$
with $\nu_n = \E(Y_{n+1} |{\cal F}_n).$

Suppose player $1$ adopts the strategy given by 
$$M_n = M = \left [ \begin{array}{r r}
 \eps  & 1-\eps \\ 
1-\eps  &  \eps  
\end{array} \right ].
$$
for some $0 < \eps < 1.$
Then $\pi_n=\pi$ with $\pi(0) = \pi(1) = 1/2$ and 
$$\theta_n = \pi_n \hat{V_n} =  -1/2$$ 
regardless of the strategy played by $2.$
 
Suppose now that  player 2 plays $Y_{n+1} = X_n$ for all $n \geq 1.$ For $\eps \neq 1/2$ hypothesis \ref{hyp:main2} is not verified and  the prediction given by (a wrong application of) theorem \ref{th:main} fails since
$$v_n \rar \sum_{x,y} \pi(x)M(x,y) U^1(x,y) = -(1-\eps).$$
\section{Proof of Theorem \ref{th:main}}
\mylabel{sec:proof}
Let $F$ denote a set--valued function mapping each point  $x \in \RR^m$ to
a  set  $F(x) \subset \RR^m.$
We call $F$ a {\em standard set valued-map} provided it verifies the three following conditions:
\bdes
\iti $F$ is a closed set-valued map. That is
$$Graph(F) = \{(x,y) \: : y \in F(x)\}$$ is a closed subset of $\RR^m
\times \RR^m.$
\itii $F$ has nonempty compact convex values, meaning that
$F(x)$ is a nonempty compact convex subset of $\RR^m$ for all $x \in \RR^m.$
\itiii
There exists $c > 0$ such that for all $x \in \RR^m$
$$\sup_{z \in F(x)} \| z\|  \leq c(1 + \| x\| )$$
where $\| \cdot \|$ denotes any  norm on $\RR^m $.
\edes
Given a standard set-valued map $F,$ set
$$F^{\delta}(u) = \{w \in \RR^m \,: \exists v \in \RR^m \, : \, d(u,v) \leq \delta,~ d(w, F(v)) \leq \delta\}.$$
The following proposition  follows from the results of Bena\"{\i}m, Hofbauer and Sorin (2005). 
\bprop
\mylabel{th:bhs}
 Let $(x_n)$ and $(U_n)$ be  discrete time processes 
living in $\mathbb{R}^m$ and $(\gamma_n)$ a sequence of nonnegative numbers.
Let $(F_n)$  be a sequence of set-valued maps and let $F$ be a standard set valued-map. 
Assume that 
\bdes
\iti 
$$
x_{n+1} - x_n  - \gamma_{n+1}  U_{n+1} \in  \gamma_{n+1} F_n(x_n)
$$
\itii
$$\sum_n \gamma_n = \infty, \: \lim_{n \to \infty} \gamma_n = 0.$$
\item 
\itiii  For all \ $T > 0$
$$\lim_{n \to \infty} \ \ \sup \left\{\left\| \sum_{i = n}^{k-1} \gamma_{i+1} U_{i+1}\right\|
   \, : k = n+1
, \ldots,  m(\tau_n + T) \right\} = 0$$
where
$$\tau_n = \sum_{i = 1}^n \gamma_i$$ and
\beq
\label{eq:tauinv} m(t) = \sup \{k \geq 0 \, : t \geq \tau_k \}
\eeq
\itiv $ \sup_n \| x_n\| = M < \infty,$
\itv For all $\delta > 0$ 
 there exists $n_0$ such that 
$$F_n(x_n) \subset F^{\delta}(x_n)$$
for all $n \geq n_0.$
\edes
  Then the limit set of $(x_n)$ is an attractor free set of the dynamics induced by $F.$
\eprop
\brem
{\em This proposition is purely deterministic. If the $(x_n), (U_n)$ are random processes, the assumptions have to be understood almost surely.
}
\erem
\brem
{\em If condition $(v)$ is strengthen to  $F_n = F,$ Proposition \ref{th:bhs} follows from  Proposition 1.3 and  Theorem 4.3 of Bena\"{\i}m, Hofbauer and Sorin (2005). Under the weaker hypothesis $(v),$ it suffices to verify that the arguments given in the proof of Proposition 1.3 adapt verbatim.}
\erem
With the notation of  the preceding sections, write
$$v_{n+1} - v_n =  \frac{1}{n+1} [-v_n + V_{n+1}] = \frac{1}{n+1} [-v_n + \theta_n + U_{n+1}]$$
where 
\beq
\label{eq:un}
U_{n+1} = V_{n+1} -\theta_n.
\eeq
Hence, conditions $(i),(ii)$ and $(iv)$ of the previous proposition are satisfied
with 
$F_n(u) = - u + C_n(u)$ and $\gamma_n = \frac{1}{n}.$
Condition $(v)$ follows from the next lemma.
\blem 
\mylabel{th:graph}
Let $C$ be adapted. For $u \in \Sigma$ and $\delta > 0$ set
$$C^{\delta}(u) = \{w \in \Sigma \,: \exists v \in \Sigma \, : \, d(u,v) \leq \delta, d(w, C(v)) \leq \delta\}.$$ Then
for all $\delta > 0$ there exists $n_0$ such that
$$C_n(u) \subset C^{\delta}(u)$$ for all $n \geq n_0$ and $u \in p(C_n).$
\elem
\prf Let $\Gamma_n = p(C_n).$ Assume to the contrary that there exist sequences $u_{n_k} \in \Gamma_{n_k}$ and $v_{n_k} \in C_{n_k}(u_{n_k})$
such that $n_k \to \infty$ and $v_{n_k} \not \in C^{\delta}(u_{n_k}).$ By compactness we may assume that $u_{n_k} \to u, v_{n_k} \to v \in C(u).$ Hence for $n_k$ large enough $d(u_{n_k},u) < \delta$ and $d(v_{n_k},v) < \delta$ proving that $v_{n_k} \in C^{\delta}(u_{n_k}).$ \qed

\medskip
To conclude the proof of theorem \ref{th:main} it remains to verify condition $(iii)$ of proposition 
\ref{th:bhs}.
\blem
Under hypothesis \ref{hyp:main} and \ref{hyp:main2}, the sequence $(U_n)$ defined by (\ref{eq:un}) verifies hypothesis $(iii)$ of proposition \ref{th:bhs}.
\elem
\prf
Set  $\frac{1}{n+1} U_{n+1} = \eps_{n+1}^0 + \eps_{n+1}$ with
$$\eps_{n+1}^0 = \frac{1}{n+1}(V_{n+1} - \hat{V}_n(X_{n+1})),$$
and 
\begin{eqnarray*}
\eps_{n+1} &=& \frac{1}{n+1}(\hat{V}_n(X_{n+1}) - \pi_n \hat{V}_n) \\
&=& \frac{1}{n+1} (Q_n - M_n Q_n)\hat{V}_n(X_{n+1})
\end{eqnarray*}where the last equality follows from the definition of $Q_n.$
Now, write $\eps_{n+1}=  \sum_{i = 1}^4 \eps^i_{n+1},$ where
$$\eps_{n+1}^1 = \frac{1}{n+1}
[ Q_n \hat{V}_n(X_{n+1}) - M_n Q_n \hat{V}_n (X_n)],$$
$$\eps_{n+1}^2 = \frac{1}{n+1} M_n Q_n \hat{V}_n(X_n) - \frac{1}{n} M_n Q_n \hat{V}_n(X_n),$$
$$\eps_{n+1}^3 = \frac{1}{n} M_n Q_n \hat{V}_n (X_n) - \frac{1}{n+1} M_{n+1} Q_{n+1} \hat{V}_{n+1}(X_{n+1}),$$
$$\eps_{n+1}^4 = \frac{1}{n+1} M_{n+1} Q_{n+1} (\hat{V}_{n+1} - \hat{V}_n)
(X_{n+1}),$$
$$\eps_{n+1}^5 = \frac{1}{n+1} [ M_{n+1} Q_{n+1} \hat{V}_n(X_{n+1}) - M_n Q_n \hat{V}_{n} (X_{n+1})].$$
For $i = 0,\ldots,5$, let
$$\eps^i_n(T) = \sup\left\{\left\|\sum_{j = n}^{k-1} \eps_{j+1}^i\right\| : k = n+1,\ldots,m(\tau_n+T)\right\}.$$
Since $\Sigma$ is compact there exists a finite constant $R$ such that $$||V_n|| + \sum_x ||\hat{V}_n(x)|| \leq R.$$

The sequence $(\eps_n^0)$ is a martingale difference with $$\E(||\eps_{n+1}^0||^2|\cF_n) \leq R^2/(n+1)^2.$$ Therefore, by Doob's convergence theorem for $L^2$ martingales, 
$\lim_{n \to \infty} \eps_n^0(T) =  0$ a.s.

The sequence $(\eps_n^1)$ is a martingale difference with
$||\eps_{n+1}^1|| \leq R |Q_n| /(n+1).$ Thus by a classical application of exponential martingale inequality (inequality (18) in Bena\"{\i}m (1999))
we have for all positive $\alpha$,
$$\P (\eps_n^1(T) \geq \alpha) \leq c 
\exp\left(\frac{-\alpha^2}{c \sum_{i = n}^{m(\tau_n+T)} (R^2 |Q_i|^2 /i^2)}\right)$$
for some positive constant $c.$
By hypothesis \ref{hyp:main}, for any $\eps > 0$ and $n$ large enough (note that $(n-1)e^T\leq m(\tau_n+T)\leq ne^T$)
$$\sum_{i = n}^{m(\tau_n+T)} (R^2 |Q_i|^2 /i^2)
\leq \sum_{i = n}^{m(\tau_n+T)}\frac{1}{i} \frac{\eps}{\log(i)}
\leq \frac{\eps(T+1)}{\log(n)}.$$
Thus 
$$\sum_{n} \P(\eps_n^1(T) \geq \alpha)) < \infty$$
and $\lim_{n \to \infty} \eps_{n+1}^1 = 0$ a.s. by Borel-Cantelli Lemma.

\smallskip
For $n+1\leq k\leq m(\tau_n+T)$,
$$\sum_{j = n}^{k-1} \eps_{j+1}^2 = \sum_{j = n}^{k-1} \frac{M_j Q_j \hat{V}_j(X_j)}{(j+1)j}.$$ Thus $$\eps^2_n(T) \leq R \sum_{j = n}^{m(\tau_n+T)} 
 \frac{|Q_j|}{j(j+1)} \leq R (T+1)  \sup_{j \geq n} \frac{|Q_j|}{j}.$$ 
By hypothesis \ref{hyp:main}, this goes to zero a.s. when $n \to \infty.$

\smallskip
For $n+1\leq k\leq m(\tau_n+T)$,
$$\sum_{j = n}^{k-1} \eps_j^3 = \frac{1}{n} M_nQ_n \hat{V}_n(X_n) - 
 \frac{1}{k}M_kQ_k \hat{V}_k(X_k),$$ so that
$$\eps^3_n(T) \leq 2R \sup_{i \geq n}\frac{1}{i} |Q_i|$$
and $\eps^3_n(T) \to 0$ a.s. as $n \to \infty$ by hypothesis \ref{hyp:main}.

\smallskip
The term $\eps^4_n(T)$ is dominated by 
$$(T+1) \sup_{i \geq n} \sup_x |M_{i+1}Q_{i+1} (\hat{V}_{i+1} - \hat{V}_i)(x)|$$
which converges a.s. towards $0$ as $n \to \infty$ by hypothesis \ref{hyp:main2}.

\smallskip
Finally, since $M_nQ_n = Q_n + \Pi_n - I$
$$\eps^5_n = \frac{1}{n+1} \left[(Q_{n+1} - Q_n) \hat{V}_n(X_{n+1}) + (\Pi_{n+1} - \Pi_n) \hat{V}_n\right].$$
Hence
$$\eps^5_n(T) \leq R  (T+1) \sup_{i \geq n} (|Q_{i+1} - Q_i| + |\pi_{i+1} - \pi_i|) \to 0$$ a.s. by hypothesis \ref{hyp:main}. This completes the proof of (iii). \qed


\begin{thebibliography}{10}
\bibitem{Ben97} M. Bena\"{\i}m, (1997), \emph{Vertex Reinforced Random Walks and a Conjecture of Pemantle},
Annals of Probability, (1997) Vol. 25, No 1, 361-392


\bibitem{ben}
M. Bena\"{\i}m, (1999), \emph{Dynamics of stochastic approximation algorithms},
 S\'{e}minaire de Probabilit\'es XXXIII, Lecture Notes in
Math. {\bf 1709}, 1--68, Springer.

\bibitem{BenHir96}
M. Bena{\"{\i}}m and M.W. Hirsch, (1996),
\emph{Asymptotic pseudotrajectories and chain recurrent flows, with
  applications}, J. Dynam. Differential Equations, 8,141--176.


\bibitem{BH99}
M. Bena\"{\i}m and M.W. Hirsch, (1999),  \emph{Mixed equilibria and dynamical
systems arising from fictitious play in perturbed games,} 
  Games and  Economic  Behavior,  {\bf 29},  36-72.
  
\bibitem{BHS}
M. Bena\"{\i}m, J. Hofbauer and S. Sorin,  (2005), \emph{Stochastic approximations and
differential inclusions},  SIAM Journal on Control and Optimization, {\bf 44}, 328-348.

\bibitem{BHS2}
M. Bena\"{\i}m, J. Hofbauer and S. Sorin,  (2006), \emph{Stochastic approximations and
differential inclusions. Part II: Applications},  Mathematics of Operations Research, 31, 673-695.

\bibitem{blr}
M. Bena\"{\i}m, M. Ledoux and O. Raimond, (2002),
\emph{Self-interacting diffusions},
 Probab. Theor. Relat. Fields \textbf{122}, 1-41.

\bibitem{br2}
M. Bena\"{\i}m and O. Raimond, (2003),
\emph{Self-interacting diffusions II: Convergence in Law.},
Annales de l'institut Henri-Poincar\'e \textbf{6}, 1043-1055. 

\bibitem{br3}
M. Bena\"{\i}m and O. Raimond, (2005),
\emph{Self-interacting diffusions III: Symmetric interactions.},
Annals of Probability {\bf 33},  no. 5, 1717--1759. 


\bibitem{BenMetPr90}
A. Benveniste and M. {M\'etivier} and P. Priouret, (1990),
``Stochastic Approximation and Adaptive Algorithms'', Springer-Verlag, Berlin and New York.

\bibitem{Bla}
D. Blackwell, (1956), \emph{An analog of the minmax theorem for vector payoffs}, Pacific Journal of Mathematics, {\bf 6}, 1-8.

\bibitem{brown}
G. Brown, (1951), \emph{Iterative solution of games by
fictitious play}, in Koopmans T.C. (ed.)  Activity Analysis of
Production and Allocation,
Wiley, 374-376.

\bibitem{Duflo96}
M. Duflo, (1996), ``Algorithmes Stochastiques'', Math\'ematiques et Applications, Springer-Verlag,
vol  23.

\bibitem{FK}
Fudenberg, D. and Kreps, K. (1993).
\emph{Learning mixed equilibria}.
 Games and Econom. Behav., 5:320--367.


\bibitem{FL}
D. Fudenberg and D. K. Levine, (1998) ``The Theory of Learning in Games'',
  MIT Press.

\bibitem{Hajek}
B. Hajek, (1982), \emph{Cooling schedules for optimal annealing}, Math. Oper. Res., 13, 2, pp 311-329.

\bibitem{hofbauer-sandholm}
J. Hofbauer and W. Sandholm., (2002) \emph{On the global convergence of stochastic fictious play}, Econometrica, 70:2265-2294.

\bibitem{holley}
R. Holley and D. Stroock, (1988),  \emph{Simulated Annealing via Sobolev Inequalities},
Commun. Math. Phys. 115, 553-568.

\bibitem{KushCla78}
H.J. Kushner and C.C. Clarck, (1978),
``Stochastic Approximation for Constrained and Unconstrained
 Systems", Springer-Verlag , Berlin and New York.

\bibitem{Ljung77}
L. Ljung, (1977), \emph{Analysis of recursive stochastic algorithms},
IEEE Trans. Automat. Control, AC-22", 551-575.

\bibitem{metpri}
M. M\'etivier and P.~ Priouret, (1987), \emph{Th\'eor\`emes de convergence presque
 sure pour une classe d'algorithmes stochastiques \`a pas d\'ecroissant},
Probability Theory and Related Fields,  74 ,403-428.

\bibitem{Miclo} 
L.~Miclo, (1992) \emph{Recuit simul\'e sans potentiel sur un ensemble fini,}  Séminaire de probabilit\'es (Strasbourg), tome 26, pp. 47-60.

\bibitem{ms}
D. Monderer and L.S. Shapley, (1996)  \emph{Fictitious play property for
  games with identical interests}, 
J.~Economic Theory,  \textbf{68}, 258--265.



\bibitem{Pem88}
R. Pemantle, (1988), ``Random processes with reinforcement'' Doctoral Dissertation, M.I.T.

\bibitem{Pem92}
R. Pemantle, (1992),
 \emph{Vertex Reinforced Random Walk},
Probab. Theor. Relat. Fields \textbf{92}, 117-136.

\bibitem{Pem07}
R. Pemantle, (2007),
\emph{A survey of random processes with reinforcement,}
Probability survey, Vol 4, 1-79.

\bibitem{Pemvol}
R. Pemantle and S. Volkov, (1999), 
\emph{Vertex-reinforced random walk on $\ZZ$ has finite range}, Annals of Probability, {\bf 27}, 1368-1388.

\bibitem{Raimond}
O. Raimond, (2008),
\emph{Self-interacting diffusions: A simulated annealing version}, To appear in Probab. Theor. Relat. Fields.

\bibitem{lsc}
L. Saloff-Coste,(1997), ``Lectures on finite Markov Chains'',
Lectures on Probability Theory and Statistics 1996, Lecture Notes in Mathematics, Vol 1665.


\bibitem{tarres}
P. Tarr\`es, (2004), 
\emph{VRRW on $\ZZ$ enventually get stuck at a set of five points}, Annals of probability 32, {\bf 3}, 1455-1478.

\bibitem{volkov}
S. Volkov, (2001),
\emph{Vertex-reinforced random walks on arbitrary graphs.} Annals of probability 29, 66-91.

\end{thebibliography}
\end{document}